# STOCHASTIC NETWORKS WITH MULTIPLE STABLE POINTS

By Nelson Antunes, Christine Fricker, Philippe Robert
and Danielle Tibi

*Universidade do Algarve, INRIA, INRIA and Université Paris 7*

This paper analyzes stochastic networks consisting of a set of finite capacity sites where different classes of individuals move according to some routing policy. The associated Markov jump processes are analyzed under a *thermodynamic limit* regime, that is, when the networks have some symmetry properties and when the number of nodes goes to infinity. An intriguing stability property is proved: under some conditions on the parameters, it is shown that, in the limit, several stable equilibrium points coexist for the empirical distribution. The key ingredient of the proof of this property is a dimension reduction achieved by the introduction of two energy functions and a convenient mapping of their local minima and saddle points. Networks with a unique equilibrium point are also presented.

**1. Introduction.** This paper studies the asymptotic behavior of a class of stochastic networks. A general description of the basic mechanisms of these systems is given in terms of a finite particle system or in terms of a queueing network.

*Description.* As for some classical stochastic processes, like the zero range process (see Liggett [17]), one can give two alternative presentations for these networks.

*A particle system.* It can be thought of as a set of sites where $K$ different types of particles coexist. At a given site, for $1 \leq k \leq K$, external particles of type $k$ with mass $A_k \in \mathbb{N}$ arrive at rate $\lambda_k$. A type $k$ particle stays an exponential time with parameter $\gamma_k$ at a site and then moves randomly to another site. A type $k$ particle leaves the system at rate $\mu_k$. Mass constraint: The total mass of particles at any site must be less than $C \in \mathbb{N}$, so that a particle arriving at a site is accepted only if this constraint is satisfied, otherwise it is rejected from the system.

*A queueing network.* It can be described as a set of identical finite capacity nodes where customers move from one node to another node uniformly chosen at random, being accepted if there is enough room and, otherwise, being rejected. If he is not rejected during his travel through the network, a customer leaves the network after his total service time. Different classes of customers access the network: Classes differ by their arrival rates, total service times, residence times at the nodes and also by the capacities they require at the nodes they visit. For example, a "light" customer will require one unit of capacity while a "heavy" customer may ask for a significant proportion of the total capacity of the node. External class $k$ customers arrive at rate $\lambda_k$ at any node. During their total service time, which ends at rate $\mu_k$, class $k$ customers are transferred from a node to another one at rate $\gamma_k$. A class $k$ customer occupies $A_k \in \mathbb{N}$ units of capacity at any visited node, if this amount of capacity is not available, he is rejected.

*Large distributed networks and statistical mechanics.* These stochastic networks have been introduced in Antunes, Fricker, Robert and Tibi [1] to represent the time evolution of a wireless network. Recent developments of mobile or sensor networks have given a strong impetus to the analysis of the associated mathematical models. See Borst, Hegde and Proutière [2], Gupta and Kumar [11] and Kermarrec, Massoulié and Ganesh [16], for example. The point of view of statistical mathematics has been introduced in the analysis of stochastic networks some time ago by Dobrushin to study various aspects of queueing networks such as departure processes of queues, capacity regions or occupancy problems. See Karpelevich, Pechersky and Suhov [12] for a survey. The networks considered quite recently have a very large number of nodes, of the order of $10^5$–$10^6$ nodes in practice, this establishes an even stronger connection with classical models of statistical physics. At the same time, due to the variety of possible topological structures and algorithms governing the behavior of these networks, new classes of mathematical models are emerging. This is clearly a promising research area for statistical mechanics methods.

*Outline of the paper.* Assuming Poisson arrivals and exponential distributions for the various service times and residence times, the time evolution of such a network is described by a Markov jump process with values in some finite (but large) state space. These associated Markov processes are, in general, not reversible, and little is known on the corresponding invariant distributions.

*Mean-field convergence.* These networks are analyzed under a thermodynamic scaling, that is, when the number of nodes $N$ goes to infinity. It is shown, Section 2, that the process of the empirical distribution $(Y^N(t))$ of the system converges to some dynamical system $(y(t))$ verifying

$$(1) \qquad \frac{d}{dt} y(t) = V(y(t)), \qquad t \geq 0,$$

where $(V(y), y \in \mathcal{P}(\mathcal{X}))$ is a vector field on $\mathcal{P}(\mathcal{X})$, the set of probability distributions on the finite set $\mathcal{X}$ defined by

$$\mathcal{X} = \{n = (n_k) \in \mathbb{N}^K : A_1 n_1 + \cdots + A_K n_K \leq C\}.$$

The dynamical system $(y(t))$ is therefore a limiting description of the original Markov process $(Y^N(t))$.

In general, the dimension of the state space $\mathcal{P}(\mathcal{X})$ of $(y(t))$ is quite large so that it is difficult to study this dynamical system in practice (an important example considered at the end of the paper is of dimension 22). The classification of the equilibrium points of $(y(t))$ with regard to the stability property is the main problem addressed in this paper. As it will be seen (cf. Proposition 5) the analysis of these points gives also insight on the limiting behavior of the invariant distribution of Markov process $(Y^N(t))$.

As a first result, it is shown that the equilibrium points are indexed by a "small" $K$-dimensional subset of $\mathbb{R}_+^K$ (the number $K$ of different classes of customers is usually quite small). They are identified as those elements of the family of probability distributions $\nu_\rho$ on $\mathcal{X}$ indexed by $\rho = (\rho_k) \in \mathbb{R}_+^K$,

$$\nu_\rho(n) = \frac{1}{Z(\rho)} \prod_{k=1}^{K} \frac{\rho_k^{n_k}}{n_k!}, \qquad n \in \mathcal{X},$$

where $Z(\rho)$ is the partition function, for which $\rho$ satisfies the fixed point equations

$$(2) \qquad \rho_k = \frac{\lambda_k + \gamma_k \sum_{m \in \mathcal{X}} m_k \nu_\rho(m)}{\gamma_k + \mu_k}, \qquad 1 \leq k \leq K.$$

The probability $\nu_\rho$ can also be seen as the invariant distribution of a multiclass $M/M/C/C$ queue. In fact the explicit limiting dynamics for the empirical distribution of the nodes is formally similar to the evolution equation for the probability distribution of the multiclass $M/M/C/C$ queue with arrival rates $\lambda_k$, service rates $\mu_k + \gamma_k$ and capacity requirements $A_k$, with the following crucial difference: the "external" arrival rates $\lambda_k$ are supplemented by "internal" arrival rates (corresponding to the mean arrival rates due to transfers from other nodes) which depend on the current state of the network.

Although the equilibrium points are indexed by a subset of $\mathbb{R}_+^K$, a dimension reduction of the dynamical system $(y(t))$ on $\mathcal{P}(\mathcal{X})$ to some dynamical system of $\mathbb{R}_+^K$ does not seem to hold. This intriguing phenomenon has also been noticed by Gibbens, Hunt and Kelly [10] for a different class of loss networks. See below.

It is shown in Section 3 that there is a unique equilibrium point when all the capacity requirements of customers are equal. For arbitrary capacity requirements, a limiting regime of the fixed point equations with respect to $\rho$ is also analyzed: The common capacity $C$ of the nodes goes to infinity and the arrival rates are proportional to $C$. In this context, Theorem 2 shows that there is essentially one unique solution: It is shown that, if $\bar{\rho}_C$ is a solution of equation (2) for capacity $C$ then, in the limit, $\bar{\rho}_C \sim \eta C$, where $\eta$ is some vector with an explicit representation in terms of the parameters of the network.

*A Lyapunov function.* In Section 4, going back to the general case, a key entropy-like function $g$ is then introduced and shown to be a Lyapunov function for $(y(t))$ so that the local minima of $g$ on $\mathcal{P}(\mathcal{X})$ correspond to the stable points of $(y(t))$. Because of the order of magnitude of $\dim \mathcal{P}(\mathcal{X})$, the identification of the local minima of $g$ on $\mathcal{P}(\mathcal{X})$ is still not simple. Using this Lyapunov function, it is proved that, if $\pi_N$ is the invariant distribution of the process of the empirical distributions $(Y^N(t))$, then the support of any limiting point of $(\pi_N)$ is carried by the set of equilibrium points of $(y(t))$. In particular, when there is a unique equilibrium point $y_\infty$ for $(y(t))$, the sequence of invariant distributions $(\pi_N)$ converges to the Dirac mass at $y_\infty$.

*A dimension reduction.* A second key function $\phi$ on the lower dimensional space $\mathbb{R}_+^K$ is introduced in Section 5. The main result of the paper, Theorem 3, establishes a one to one correspondence between the local minima of $g$ on $\mathcal{P}(\mathcal{X})$ and the local minima of $\phi$. The dimension reduction for the problem of stability of the equilibrium points is therefore achieved not through dynamical systems but through the energy functions $g$ and $\phi$. This result is interesting in its own right and seems to be a promising way of studying other classes of large networks.

*Phenomenon of bistability for $(y(t))$.* With these results, an example of a network with two classes of customers and at least *three zeroes* for $V$ is exhibited in Section 6, two of them being "stable" and the other one being a saddle point. In this case the asymptotic dynamical system $(y(t))$ has therefore a bistability property. This suggests the following (conjectured) bistability property for the original process describing the state of the network: it lives for a very long time in a region corresponding to one of the stable points and then, due to some rare events, it reaches, via a saddle point, the region of another stable point and so on. This conjectured phenomenon is known as metastability in statistical physics. A formal proof of this phenomenon seems to be quite difficult to obtain. The only tools available in this domain use either the Gibbsian characteristics of the dynamics (see Olivieri and Vares [19]) or at least the reversibility of the Markov process, Bovier [3, 4]. None of these properties holds in our case.

Note that in Antunes et al. [1] it is proved that, for similar networks under Kelly's scaling, there is a unique equilibrium point. Contrary to the model considered here, the capacity requirement of a customer in [1] does not depend on his class. On the other hand, the networks analyzed here have a symmetrical structure: all the nodes have the same capacity and the routing is uniform among all the other nodes. In Antunes et al. [1], the routing mechanisms are quite general.

*Multiple stable equilibrium points and local dynamics in stochastic networks.* Asymptotic dynamical systems with multiple stable points are quite rare in stochastic networks. Gibbens, Hunt and Kelly [10] have shown, via an approximated model, that such an interesting phenomenon may occur in a loss network with a rerouting policy. Marbukh [18] analyzes, also through some approximation, the bistability properties of similar loss networks. The dynamics of the networks of Gibbens, Hunt and Kelly [10], Marbukh [18] and of this paper are local in the following sense: The interaction between two nodes only depends on the states of these two nodes and not on the state of the whole network.

A key feature of these models is the subtle interplay between the local description of the dynamics and its impact on the macroscopic state of the network through the existence of several equilibrium points. In statistical physics these phenomena have been known for some time. See den Hollander [8], Olivieri and Vares [19] and the references therein for a general presentation of these questions. See also Bovier [3, 4] for a potential theoretical approach in the case of reversible Markov processes and Catoni and Cerf [6] for a study of the saddle points of perturbed Markov chains. For the more classical setting of global dynamics, the large deviation approach is developed in Freidlin and Wentzell [9].

*Phase transitions in uncontrolled loss networks.* This is a related topic. It is known that, for some loss networks on infinite graphs, there may be several equilibrium distributions. If the loss network is restricted to a finite subgraph $\mathcal{F}$ of the infinite graph $\mathcal{G}$, its equilibrium distribution $\pi_\mathcal{F}$ is uniquely determined. It turns out that, depending on the boundary conditions on $\mathcal{F}$, the sequence $(\pi_\mathcal{F}, \mathcal{F} \subset \mathcal{G})$, may have distinct limiting values which are invariant distributions for the case of the infinite graph. Significant results have already been obtained in this domain. See Spitzer [22], Kelly [15] for a survey, Zachary [24], Zachary and Ziedins [25] and Ramanan, Sengupta, Ziedins and Mitra [20].

**2. The asymptotic dynamical system.** Two nodes $i, j \in \{1, \ldots, N\}$ of the network interact through the exchange of customers at rate of the order of $1/N$. Due to the symmetrical structure of the network, a stronger statement holds: the impact on $i$ of all nodes different from $i$ appears only through some averaged quantity. For $1 \leq k \leq K$, the input rate of class $k$ customers at node $i$ from the other nodes is

$$\frac{1}{N-1} \sum_{1 \leq j \leq N, j \neq i} \gamma_k X_{j,k}^N(t).$$

If this quantity is close to $\gamma_k \mathbb{E}(X_{1,k}^N(t))$, a *mean field property* is said to hold. Note that, if the network starts from some symmetrical initial state, the random variables $X_{j,k}^N(t)$, $j = 1, \ldots, N$, have the same distribution.

THEOREM 1. *If $Y^N(0)$ converges weakly to $z \in \mathcal{P}(\mathcal{X})$ as $N$ tends to infinity, then $(Y^N(t))$ converges in the Skorohod topology to the solution $(y(t))$ of the ordinary differential equation*

$$(3) \qquad y'(t) = V(y(t)),$$

*where $(y(t))$ is the solution starting from $y(0) = z$ and, for $y \in \mathcal{P}(\mathcal{X})$, the vector field $V(y) = (V_n(y), n \in \mathcal{X})$ on $\mathcal{P}(\mathcal{X})$ is defined by*

$$
\begin{aligned}
V_n(y) &= \sum_{k=1}^{K} (\lambda_k + \gamma_k \langle \mathbb{I}_k, y \rangle)(y_{n-f_k} \mathbb{1}_{\{n_k \geq 1\}} - y_n \mathbb{1}_{\{n+f_k \in \mathcal{X}\}}) \\
&\quad + \sum_{k=1}^{K} (\gamma_k + \mu_k)((n_k + 1) y_{n+f_k} \mathbb{1}_{\{n+f_k \in \mathcal{X}\}} - n_k y_n),
\end{aligned}
$$
(4)

*where $\langle \mathbb{I}_k, y \rangle = \sum_{m \in \mathcal{X}} m_k y_m$ and $f_k$ is the $k$th unit vector of $\mathbb{R}^K$.*

By convergence in the Skorohod topology, one means the convergence in distribution for Skorohod topology on the space of trajectories.

Note that equation (4) gives the derivative $dy_n(t)/dt = V_n(y(t))$ of $y_n(t)$ as increasing proportionally to the difference $y_{n-f_k} - y_n$ by some factor $\lambda_k + \gamma_k \langle \mathbb{I}_k, y \rangle$ which measures the speed at which nodes in state $n - f_k$ turn to state $n$ (due to an arrival of some type $k$ customer). In this factor, $\gamma_k \langle \mathbb{I}_k, y \rangle$ is added to the external arrival rate $\lambda_k$ of class $k$ customers at any node, and hence appears as the internal arrival rate of class $k$ customers at any node. This feature characterizes the mean field property. Indeed, $\langle \mathbb{I}_k, y \rangle$ is the mean number of class $k$ customers per node when the empirical distribution of the $N$ nodes is $y$; so that $\gamma_k \langle \mathbb{I}_k, y \rangle$ is the mean emission rate per node of class $k$ customers to the rest of the network.

PROOF OF THEOREM 1. The martingale characterization of the Markov jump process $(Y_n^N(t))$, see Rogers and Williams [21], shows that

$$M_n^N(t) = Y_n^N(t) - Y_n^N(0) - \int_0^t \sum_{w \in \mathcal{P}(\mathcal{X}) \setminus \{Y^N(s)\}} \Omega_N(Y^N(s), w)(w_n - Y_n^N(s))\, ds$$

is a martingale with respect to the natural filtration associated to the Poisson processes involved in the arrival processes, service times and residence times. By using the explicit expression of the $Q$-matrix $\Omega_N$, trite (and careful) calculations finally show that the following relation holds:

$$
\begin{aligned}
Y_n^N(t) = {}& Y_n^N(0) + M_n^N(t) \\
& + \int_0^t \sum_{k=1}^K \left( \lambda_k + \frac{\gamma_k N}{N-1} \sum_{m \in \mathcal{X}} m_k Y_m^N(s) \right) \\
& \qquad \times (Y_{n-f_k}^N(s) \mathbb{1}_{\{n_k \geq 1\}} - Y_n^N(s) \mathbb{1}_{\{n+f_k \in \mathcal{X}\}})\, ds \\
& + \int_0^t \sum_{k=1}^K (\gamma_k + \mu_k)((n_k+1) Y_{n+f_k}^N(s) \mathbb{1}_{\{n+f_k \in \mathcal{X}\}} - n_k Y_n^N(s))\, ds \\
& + \int_0^t \sum_{k=1}^K \frac{\gamma_k}{N-1}(n_k Y_n^N(s) - (n_k-1) Y_{n-f_k}^N(s) \mathbb{1}_{\{n_k \geq 1\}})\, ds.
\end{aligned}
$$
(5)

From there, with a similar method as in Darling and Norris [7], it is not difficult to prove that if $Y^N(0)$ converges to $z$ then:

— the sequence $(Y^N(t))$ of process is tight for the Skorohod topology;
— any limit $(y(t))$ of $(Y^N(t))$ is continuous and satisfies the following deterministic differential equation, $y(0) = z$ and

$$
\begin{aligned}
y_n'(t) = {}& \sum_{k=1}^K \left[ \lambda_k + \gamma_k \sum_{m \in \mathcal{X}} m_k y_m(t) \right] [y_{n-f_k}(t) \mathbb{1}_{\{n_k \geq 1\}} - y_n(t) \mathbb{1}_{\{n+f_k \in \mathcal{X}\}}] \\
& + \sum_{k=1}^K (\gamma_k + \mu_k)[(n_k+1) y_{n+f_k}(t) \mathbb{1}_{\{n+f_k \in \mathcal{X}\}} - n_k y_n(t)].
\end{aligned}
$$

This is exactly equation (4). The uniqueness of the solution of this differential equation implies that such a limiting point $(y(t))$ is necessarily unique and therefore that $(Y^N(t))$ converges in distribution to $(y(t))$. The proposition is proved. □

The equilibrium points of the dynamical system defined by equation (4) are the probability distributions $y \in \mathcal{P}(\mathcal{X})$ on $\mathcal{X}$ such that $V(y)$ is zero. This condition can be written as follows: For $n \in \mathcal{X}$,

$$
\left(\sum_{k=1}^{K}[(\lambda_k + \gamma_k \langle \mathbb{I}_k, y \rangle)\mathbb{1}_{\{n+f_k \in \mathcal{X}\}} + (\gamma_k + \mu_k)n_k]\right)y_n
$$

(6)
$$
= \sum_{k=1}^{K}(\lambda_k + \gamma_k \langle \mathbb{I}_k, y \rangle)y_{n-f_k}\mathbb{1}_{\{n_k \geq 1\}}
$$
$$
+ (\gamma_k + \mu_k)(n_k+1)y_{n+f_k}\mathbb{1}_{\{n+f_k \in \mathcal{X}\}}.
$$

These equations are equivalent to local balance equations for the numbers of customers of a classical $M/M/C/C$ queue with $K$ classes of customers such that, for $1 \leq k \leq K$, class $k$ customers:

— arrive at rate $\lambda_k + \gamma_k \langle \mathbb{I}_k, y \rangle$;
— are served at rate $\gamma_k + \mu_k$;
— require capacity $A_k$.

Consequently, $(y_n)$ is the invariant distribution of this queue. It is well known (see Kelly [13], e.g.) that necessarily

(7)
$$
y_n = \nu_\rho(n) \stackrel{\text{def.}}{=} \frac{1}{Z(\rho)}\prod_{k=1}^{K}\frac{\rho_k^{n_k}}{n_k!}, \qquad n \in \mathcal{X},
$$

where, for $1 \leq k \leq K$, $\rho_k$ is the ratio of the $k$th arrival and service rates,

(8)
$$
\rho_k = \frac{\lambda_k + \gamma_k \langle \mathbb{I}_k, \nu_\rho \rangle}{\gamma_k + \mu_k},
$$

where $\langle \mathbb{I}_k, \nu_\rho \rangle$ is the average value of the $k$th component under the probability distribution $\nu_\rho$ on $\mathcal{X}$ and $Z(\rho)$ is the normalization constant; in other words, the partition function

$$
Z(\rho) = \sum_{n \in \mathcal{X}}\prod_{k=1}^{K}\frac{\rho_k^{n_k}}{n_k!}.
$$

The equilibrium points of the dynamical system $(y(t))$ are indexed by $\mathbb{R}_+^K$ whose dimension is much smaller than $\mathcal{P}(\mathcal{X})$ thereby suggesting a possible simpler description of the asymptotic behavior of the network. Despite this quite appealing perspective, it turns out that such a dimension reduction cannot be achieved directly since the subset $\{\nu_\rho, \rho \in \mathbb{R}_+^K\}$ of $\mathcal{P}(\mathcal{X})$ is not left invariant by the dynamical system $(y(t))$.

Denote by $B_k(\rho)$ the blocking probability of a class $k$ customer in an $M/M/C/C$ queue at equilibrium with $K$ classes and loads $\rho_1, \ldots, \rho_K$, that is,

$$
B_k(\rho) = \frac{1}{Z(\rho)}\sum_{n:n+f_k \notin \mathcal{X}}\prod_{h=1}^{K}\frac{\rho_h^{n_h}}{n_h!},
$$

it is easily checked that $\langle \mathbb{I}_k, \nu_\rho \rangle = \rho_k(1 - B_k(\rho))$, so equation (8) becomes then

$$\mu_k \rho_k = \lambda_k - \gamma_k \rho_k B_k(\rho).$$

The following proposition summarizes these results.

PROPOSITION 1. *The equilibrium points of the dynamical system $(y(t))$ defined by equation (3) are exactly the probability distributions $\nu_\rho$ on $\mathcal{X}$,*

$$(9) \qquad \nu_\rho(n) = \frac{1}{Z(\rho)} \prod_{\ell=1}^{K} \frac{\rho_\ell^{n_\ell}}{n_\ell!}, \qquad n = (n_\ell) \in \mathcal{X},$$

*where $Z(\rho)$ is the partition function,*

$$Z(\rho) = \sum_{m \in \mathcal{X}} \prod_{\ell=1}^{K} \frac{\rho_\ell^{m_\ell}}{m_\ell!}$$

*and $\rho = (\rho_k, 1 \leq k \leq K)$ is a vector of $\mathbb{R}_+^K$ satisfying the system of equations*

$$(10) \qquad \lambda_k = \rho_k \left( \mu_k + \gamma_k \sum_{n\,:\,n+f_k \notin \mathcal{X}} \prod_{\ell=1}^{K} \frac{\rho_\ell^{n_\ell}}{n_\ell!} \bigg/ \sum_{n \in \mathcal{X}} \prod_{\ell=1}^{K} \frac{\rho_\ell^{n_\ell}}{n_\ell!} \right), \qquad 1 \leq k \leq K.$$

*There always exists at least one equilibrium point.*

PROOF. Only the existence result has to be proved. According to equations (7) and (8), $y$ is an equilibrium point if and only if it is a fixed point of the function

$$\mathcal{P}(\mathcal{X}) \longrightarrow \mathcal{P}(\mathcal{X}),$$

$$y \longrightarrow \nu_{\rho(y)},$$

with $\rho(y) = (\rho_k(y))$ and, for $1 \leq k \leq K$, $\rho_k(y) = (\lambda_k + \gamma_k \langle \mathbb{I}_k, y \rangle)/(\mu_k + \gamma_k)$. This functional being continuous on the convex compact set $\mathcal{P}(\mathcal{X})$, it necessarily has a fixed point by Brouwer's theorem. □

A similar situation occurs in Gibbens, Hunt and Kelly [10] where the equilibrium points are also indexed by the solutions $\rho \in \mathbb{R}_+$ of some fixed point equations and the bistability properties of the system are analyzed through numerical estimates. Here, a detailed study of the stability properties of the equilibrium points is achieved.

*Notation.* In the following, we will denote

$$\frac{\rho^n}{n!} = \prod_{k=1}^{K} \frac{\rho_k^{n_k}}{n_k!},$$

for $n = (n_k) \in \mathcal{X}$ and $\rho = (\rho_k) \in \mathbb{R}_+^K$. The system of equations (10) can then be rewritten as

$$\lambda_k = \rho_k \left( \mu_k + \gamma_k \frac{\sum_{n \,:\, n+f_k \notin \mathcal{X}} \rho^n/n!}{\sum_{n \in \mathcal{X}} \rho^n/n!} \right), \qquad 1 \leq k \leq K.$$

**3. Uniqueness results.** In this section, several situations in which the asymptotic dynamical system (3) has a unique equilibrium point, that is, when the fixed point equation (10) have a unique solution, are presented.

3.1. *Networks with constant requirements.* It is assumed that all classes of customers require the same capacity, that is, that $A_k = A$ for $k = 1, \ldots, K$. By replacing $C$ by $\lfloor C/A \rfloor$, it can be assumed that $A = 1$. In this case, if $|n|$ denotes the sum of the coordinates of $n \in \mathcal{X}$,

$$\begin{aligned}
B_k(\rho) &= \sum_{n \,:\, n+f_k \notin \mathcal{X}} \frac{\rho^n}{n!} \Big/ \sum_{n \in \mathcal{X}} \frac{\rho^n}{n!} \\
&= \sum_{|n|=C} \frac{\rho^n}{n!} \Big/ \sum_{n \in \mathcal{X}} \frac{\rho^n}{n!} \\
&= \frac{1}{C!} \left( \sum_{k=1}^{K} \rho_k \right)^C \Big/ \sum_{\ell=0}^{C} \frac{1}{\ell!} \left( \sum_{k=1}^{K} \rho_k \right)^\ell \\
&\stackrel{\text{def.}}{=} B_1\left( \sum_k \rho_k \right),
\end{aligned}$$

$B_1(\theta)$ can be represented as the stationary blocking probability of an $M/M/C/C$ queue with one class of customers and arrival rate $\theta$ and service rate 1.

In this case, fixed point equation (10) are

(11) $\qquad \rho_k = \lambda_k/(\mu_k + \gamma_k B_1(S)), \qquad k = 1, \ldots, K,$

with $S = \rho_1 + \cdots + \rho_K$. By summing up these equations, one gets that $S$ is the solution of the equation

$$S = \sum_{k=1}^{K} \frac{\lambda_k}{\mu_k + \gamma_k B_1(S)}.$$

It is easily checked that $S \to B_1(S)$ is nondecreasing and therefore that the above equation has a unique solution. The uniqueness of the vector $(\rho_k)$ follows from equations (11). The following proposition has been proved.

PROPOSITION 2. *The asymptotic dynamical system $(y(t))$ has a unique equilibrium point when capacity requirements are equal.*

In particular, when there is only one class of customers, there is a unique solution to equation (10).

3.2. *A limiting regime of fixed point equations.* Here, the fixed point equations (10) are analyzed under Kelly's scaling, that is, when the capacity $C$ goes to infinity and the arrival rates are proportional to $C$, of the order of $\lambda_k C$ for the $k$th class. For $1 \leq k \leq K$, the total service rate $\mu_k$ and the rate of residence time $\gamma_k$ are kept fixed. It will be shown that, in this case, there is a unique equilibrium point. Let $(\mathcal{N}_k, 1 \leq k \leq K)$ be a sequence of $K$ independent Poisson processes with intensity 1. As usual $\mathcal{N}_k(A)$ will denote the number of points of the $k$th process in the subset $A$ of $\mathbb{R}_+$. For $C > 0$, denote by $\rho_C = (\rho_k(C), 1 \leq k \leq K)$, a solution of the fixed point equations

$$\lambda_k C = \rho_k(C)\left(\mu_k + \gamma_k \sum_{n:\,n+f_k \notin \mathcal{X}} \frac{\rho_C^n}{n!} \bigg/ \sum_{n \in \mathcal{X}} \frac{\rho_C^n}{n!}\right), \qquad 1 \leq k \leq K,$$

this can be rewritten as

$$(12) \quad \lambda_k = \frac{\rho_k(C)}{C}\left(\mu_k + \gamma_k - \gamma_k \frac{\mathbb{P}(C - \sum_{i=1}^K A_i \mathcal{N}_i([0, \rho_i(C)]) \geq A_k)}{\mathbb{P}(C - \sum_{i=1}^K A_i \mathcal{N}_i([0, \rho_i(C)]) \geq 0)}\right).$$

This problem is related to the limit of the loss probabilities investigated and solved by Kelly [14] in a general setting in terms of an optimization problem. The proposition below is a consequence of Kelly's result.

PROPOSITION 3 (Kelly's scaling). *If $(\rho_k(C), 1 \leq k \leq K) \in \mathbb{R}_+^K$ is such that*

$$\lim_{C \to +\infty} \rho_k(C)/C = \bar{\rho}_k, \qquad 1 \leq k \leq K,$$

*with $(\bar{\rho}_k) \in \mathbb{R}_+^K$ and*

$$(13) \qquad \bar{\rho}_1 A_1 + \bar{\rho}_2 A_2 + \cdots + \bar{\rho}_K A_K \geq 1,$$

*then, for $a \in \mathbb{N}$,*

$$\lim_{C \to +\infty} \frac{\mathbb{P}(C - \sum_{k=1}^K A_k \mathcal{N}_k([0, \rho_k(C)]) \geq a)}{\mathbb{P}(C - \sum_{k=1}^K A_k \mathcal{N}_k([0, \rho_k(C)]) \geq 0)} = e^{-\omega a},$$

*where $\omega$ is the unique nonnegative solution of the equation*

$$(14) \qquad \bar{\rho}_1 A_1 e^{-\omega A_1} + \bar{\rho}_2 A_2 e^{-\omega A_2} + \cdots + \bar{\rho}_K A_K e^{-\omega A_K} = 1.$$

The main result of this section can now be stated. Basically, it states that, under Kelly's scaling, the fixed point equation (10) have a unique solution when the capacity gets large.

THEOREM 2. *If $\mu_k > 0$ for all $1 \leq k \leq K$ and if for any $C > 0$ the vector $(\rho_k(C, \lambda C))$ is any solution of equation (10) then, for $1 \leq k \leq K$,*

$$\lim_{C \to +\infty} \frac{\rho_k(C, \lambda C)}{C} = \frac{\lambda_k}{\mu_k + \gamma_k - \gamma_k e^{-\omega A_k}},$$

*where $\omega \geq 0$ is defined as*

$$\omega = \inf\left\{x \geq 0 : \sum_{k=1}^{K} \frac{\lambda_k A_k e^{-x A_k}}{\mu_k + \gamma_k - \gamma_k e^{-x A_k}} \leq 1\right\}.$$

PROOF. For $1 \leq k \leq K$, the function $C \to \bar{\rho}_k(C) = \rho_k(C, \lambda C)/C$ is bounded by $\lambda_k/\mu_k$. By taking a subsequence, it can be assumed that $\bar{\rho}_k(C)$ converges to some finite $\bar{\rho}_k$ as $C$ goes to infinity. Under the condition

$$A_1 \frac{\lambda_1}{\mu_1} + A_2 \frac{\lambda_2}{\mu_2} + \cdots + A_K \frac{\lambda_K}{\mu_K} \geq 1,$$

then necessarily $\bar{\rho}_1 A_1 + \bar{\rho}_2 A_2 + \cdots + \bar{\rho}_K A_K \geq 1$, otherwise one would have, via the law of large numbers for Poisson processes, for $a \geq 0$,

$$(15) \qquad \lim_{C \to +\infty} \mathbb{P}\left(C - \sum_{i=1}^{K} A_i \mathcal{N}_i([0, C\bar{\rho}_i(C)]) \geq a\right) = 1,$$

and equation (12) would then give the relation $\bar{\rho}_k = \lambda_k/\mu_k$ for $1 \leq k \leq K$, so that

$$\bar{\rho}_1 A_1 + \bar{\rho}_2 A_2 + \cdots + \bar{\rho}_K A_K \geq 1,$$

contradiction. From Proposition 3 and equation (12), one gets that

$$\lambda_k = \bar{\rho}_k(\mu_k + \gamma_k - \gamma_k e^{-\omega A_k}), \qquad 1 \leq k \leq K,$$

were $\omega$ is the solution of equation (14) associated to $(\bar{\rho}_k)$. Equation (14) can then be rewritten as

$$\sum_{1}^{K} \frac{\lambda_k A_k e^{-\omega A_k}}{\mu_k + \gamma_k - \gamma_k e^{-\omega A_k}} = 1.$$

The statement of the theorem is proved in this case.

Now, if it is assumed that $A_1 \lambda_1/\mu_1 + A_2 \lambda_2/\mu_2 + \cdots + A_K \lambda_K/\mu_K < 1$ then, since $\bar{\rho}_k \leq \lambda_k/\mu_k$ for all $k$, relation (15) holds and equation (12) finally gives that $\bar{\rho}_k = \lambda_k/\mu_k$, $1 \leq k \leq K$. The theorem is proved. □

**4. An energy function on $\mathcal{P}(\mathcal{X})$.** In this section, a Lyapunov function is introduced. As it will be seen, it plays a key role in the analysis of the fixed points of the asymptotic dynamical system. Define the function $g$ on the set $\mathcal{P}(\mathcal{X})$ of probability distributions on $\mathcal{X}$,

$$(16) \quad g(y) = \sum_{n \in \mathcal{X}} y_n \log(n! y_n) - \sum_{k=1}^{K} \int_0^{\langle \mathbb{I}_k, y \rangle} \log \frac{\lambda_k + \gamma_k x}{\mu_k + \gamma_k} \, dx, \qquad y \in \mathcal{P}(\mathcal{X}).$$

Recall that, for $y \in \mathcal{P}(\mathcal{X})$, $\langle \mathbb{I}_k, y \rangle = \sum_{m \in \mathcal{X}} m_k y_m$ and $\overset{\circ}{\mathcal{P}}(\mathcal{X})$ denotes the interior of the set $\mathcal{P}(\mathcal{X})$.

PROPOSITION 4. *The function $g$ is* a Lyapunov function *for the asymptotic dynamical system $(y(t))$, that is,*

$$\langle V(y), \nabla g(y) \rangle = \sum_{n \in \mathcal{X}} V_n(y) \frac{\partial g}{\partial y_n}(y) \leq 0 \qquad \forall y \in \overset{\circ}{\mathcal{P}}(\mathcal{X}),$$

*and, for $y \in \overset{\circ}{\mathcal{P}}(\mathcal{X})$, the following assertions are equivalent:*

(a) $\langle V(y), \nabla g(y) \rangle = 0$;
(b) *The coordinates of $\nabla g(y)$ are equal;*
(c) *$y$ is an equilibrium point of $(y(t))$, that is, $V(y) = 0$.*

PROOF. The vector field $V(y) = (V_n(y))$ can be written as

$$V_n(y) = \sum_{k=1}^{K} [(\lambda_k + \gamma_k \langle \mathbb{I}_k, y \rangle) y_{n-f_k} \mathbb{1}_{\{n_k \geq 1\}} + (\mu_k + \gamma_k)(n_k + 1) y_{n+f_k} \mathbb{1}_{\{n+f_k \in \mathcal{X}\}}$$
$$- ((\lambda_k + \gamma_k \langle \mathbb{I}_k, y \rangle) \mathbb{1}_{\{n+f_k \in \mathcal{X}\}} + (\mu_k + \gamma_k) n_k) y_n]$$
$$= \sum_{k=1}^{K} (F_{n+f_k}^k(y) - F_n^k(y))$$

where, for $n \in \mathcal{X}$, $F_n^k(y) = (\mu_k + \gamma_k) n_k y_n - (\lambda_k + \gamma_k \langle \mathbb{I}_k, y \rangle) y_{n-f_k} \mathbb{1}_{\{n_k \geq 1\}}$ and $F_n^k(y) = 0$ when $n \notin \mathcal{X}$; note that $F_n^k = 0$ whenever $n_k = 0$.

For $y \in \overset{\circ}{\mathcal{P}}(\mathcal{X})$,

$$\langle V(y), \nabla g(y) \rangle = \sum_{n \in \mathcal{X}} \sum_{k=1}^{K} \frac{\partial g}{\partial y_n}(y)(F_{n+f_k}^k(y) - F_n^k(y))$$
$$= \sum_{k=1}^{K} \sum_{n \in \mathcal{X}} F_n^k(y) \left( \frac{\partial g}{\partial y_{n-f_k}}(y) - \frac{\partial g}{\partial y_n}(y) \right).$$

Since, for $n \in \mathcal{X}$,

$$\frac{\partial g}{\partial y_n}(y) = 1 + \log(n! y_n) - \sum_{k=1}^{K} n_k \log \frac{\lambda_k + \gamma_k \langle \mathbb{I}_k, y \rangle}{\mu_k + \gamma_k},$$

one finally gets that the relation

$$(17) \qquad \langle V(y), \nabla g(y) \rangle = \sum_{k=1}^{K} \sum_{n \in \mathcal{X}} F_n^k(y) \log \frac{(\lambda_k + \gamma_k \langle \mathbb{I}_k, y \rangle) y_{n-f_k}}{(\mu_k + \gamma_k) n_k y_n}$$

holds. The quantity $\langle V(y), \nabla g(y) \rangle$ is thus clearly nonpositive. On the other hand, for $k$ and $n$ such that $n_k \geq 1$,

$$\frac{\partial g}{\partial y_n}(y) - \frac{\partial g}{\partial y_{n-f_k}}(y) = \log\left(\frac{(\mu_k + \gamma_k) n_k y_n}{(\lambda_k + \gamma_k \langle \mathbb{I}_k, y \rangle) y_{n-f_k}}\right),$$

hence $\langle V(y), \nabla g(y) \rangle$ is zero if and only if the coordinates of $\nabla g(y)$ are equal and this is equivalent to the system of equations

$$(\lambda_k + \gamma_k \langle \mathbb{I}_k, y \rangle) y_{n-f_k} = (\mu_k + \gamma_k) n_k y_n$$

for all $k$ and $n$ such that $n_k \geq 1$, so that $y$ is an equilibrium point of the asymptotic dynamical system. The equivalence (a), (b) and (c) is proved. □

*Convergence of the stationary distribution.* If $F$ is some real-valued function on $\mathbb{R}^{\mathcal{X}}$ and $y \in \mathcal{P}(\mathcal{X})$, the functional operator associated to the $Q$-matrix $\Omega_N$ is given by

$$\Omega_N(F)(y) = \sum_{z \in \mathcal{P}(\mathcal{X}) \setminus \{y\}} \Omega_N(y, z)(F(z) - F(y))$$

$$= \sum_{n \in \mathcal{X}} \left[ \sum_{k=1}^{K} \lambda_k y_n N \mathbb{1}_{\{n+f_k \in \mathcal{X}\}} \left( F\left(y + \frac{1}{N}(e_{n+f_k} - e_n)\right) - F(y) \right) \right.$$

$$+ \sum_{k=1}^{K} \mu_k n_k y_n N \left( F\left(y + \frac{1}{N}(e_{n-f_k} - e_n)\right) - F(y) \right)$$

$$+ \sum_{\substack{1 \leq k \leq K \\ m \in \mathcal{X}}} \frac{\gamma_k N}{N-1} n_k y_n (N y_m - \mathbb{1}_{\{m=n\}})$$

$$\times \left( F\left(y + \frac{e_{n-f_k} - e_n}{N} + \frac{e_{m+f_k} - e_m}{N} \mathbb{1}_{\{m+f_k \in \mathcal{X}\}}\right) \right.$$

$$\left.\left. - F(y) \right) \right].$$

If it is assumed that the function $F$ is of class $C^2$ on $\mathbb{R}^K$, then it is easy to check that the sequence $(\Omega_N(F)(y))$ converges to the following expression:

$$\sum_{n\in\mathcal{X}}\left[\sum_{k=1}^{K}\lambda_k y_n \mathbb{1}_{\{n+f_k\in\mathcal{X}\}}\langle\nabla F(y),e_{n+f_k}-e_n\rangle + \sum_{k=1}^{K}\mu_k n_k y_n\langle\nabla F(y),e_{n-f_k}-e_n\rangle\right.$$

$$+\sum_{\substack{1\leq k\leq K \\ m\in\mathcal{X}}}\gamma_k n_k y_n y_m(\langle\nabla F(y),e_{n-f_k}-e_n\rangle$$

$$\left.+\langle\nabla F(y),e_{m+f_k}-e_m\rangle\mathbb{1}_{\{m+f_k\in\mathcal{X}\}})\right]$$

which is defined as $\Omega_\infty(F)(y)$. Moreover, by using Taylor's formula at the second order, this convergence is *uniform with respect to* $y\in\mathcal{P}(\mathcal{X})$. By Theorem 1, one necessarily has

$$\Omega_\infty(F)(y) = \langle\nabla F(y), V(y)\rangle, \qquad y\in\mathcal{P}(\mathcal{X}).$$

Note that this identity can also be checked directly with the above equation.

PROPOSITION 5. *If $\pi_N$ denotes the invariant probability distribution of $(Y^N(t))$ on $\mathcal{P}(\mathcal{X})$, then any limiting point of $(\pi_N)$ is a probability distribution carried by the equilibrium points of the asymptotic dynamical system $(y(t))$ of Theorem 1.*

*In particular, if $(y(t))$ has a unique equilibrium point $y_\infty$, then the sequence of invariant distributions $(\pi_N)$ converges to the Dirac mass at $y_\infty$.*

PROOF. The set $\mathcal{P}(\mathcal{X})$ being compact, the sequence of distributions $(\pi_N)$ is relatively compact. Let $\widetilde{\pi}$ be the limit of some subsequence $(\pi_{N_p})$. If $F$ is a function of class $C^2$ on $\mathbb{R}^\mathcal{X}$, then for $p\geq 0$,

$$\int_{\mathcal{P}(\mathcal{X})}\Omega_{N_p}(F)(y)\pi_{N_p}(dy) = 0.$$

The uniform convergence of $\Omega_{N_p}(F)$ to $\Omega_\infty(F)$ implies that

$$0 = \int_{\mathcal{P}(\mathcal{X})}\Omega_\infty(F)(y)\widetilde{\pi}(dy) = \int_{\mathcal{P}(\mathcal{X})}\langle\nabla F(y),V(y)\rangle\widetilde{\pi}(dy),$$

so that $\widetilde{\pi}$ is an invariant distribution of the (deterministic) Markov process associated to the infinitesimal generator $\Omega_\infty$.

For $t\geq 0$, denote (temporarily) by $(y(x,t))$ the dynamical system starting from $y(0)=x\in\mathcal{P}(\mathcal{X})$. Assume that there exist $x\in\mathcal{P}(\mathcal{X})$ and $s>0$ such that $y(x,s)\in\partial\mathcal{P}(\mathcal{X})$, that is, there exists $n\in\mathcal{X}$ such that $y_n(x,s)=0$. Since $(y_n(x,t))$ is nonnegative and since the function $t\to y(x,t)$ is of class $C^1$, it implies that $V_n(y(x,s)) = \dot{y}_n(x,s) = 0$. This last relation, (4) defining the vector field $(V_n(y))$ and the equation $y_n(x,s) = 0$ give that $y_{n\pm f_k}(x,s) = 0$ for any $k$ such that $n \pm f_k \in \mathcal{X}$ and consequently, by repeating the argument, all the coordinates of $y(x,s)$ are null. Contradiction since $y(x,s)$ is a probability distribution on $\mathcal{X}$. Hence, the boundary $\partial \mathcal{P}(\mathcal{X})$ of $\mathcal{P}(\mathcal{X})$ cannot be reached in positive time by $(y(x,t))$. This entails, in particular, that $\partial \mathcal{P}(\mathcal{X})$ is negligible for any invariant distribution of $(y(x,t))$.

For $x \in \mathcal{P}(\mathcal{X})$ and $0 < s' < s$, since the function $g(y(x,\cdot))$ is of class $C^1$ on $[s', s]$ and its derivative is $\langle \nabla(g)(y(x,\cdot)), V(y(x,\cdot)) \rangle$, one has

$$(18) \qquad g(y(x,s)) - g(y(x,s')) = \int_{s'}^{s} \Omega_\infty(g)(y(x,u))\, du.$$

By integrating with respect to $\widetilde{\pi}$ this relation, the invariance of $\widetilde{\pi}$ for the process $(y(x,t))$ and Fubini's theorem show that

$$\int_{\mathcal{P}(\mathcal{X})} g(y(x,s))\widetilde{\pi}(dx) - \int_{\mathcal{P}(\mathcal{X})} g(y(x,s'))\widetilde{\pi}(dx)$$
$$= 0 = \int_{\mathcal{P}(\mathcal{X})} \int_{s'}^{s} \Omega_\infty(g)(y(x,u))\, du\, \widetilde{\pi}(dx)$$
$$= (s-s') \int_{\mathcal{P}(\mathcal{X})} \Omega_\infty(g)(x)\widetilde{\pi}(dx).$$

The integrand $\Omega_\infty(g)(x) = \langle \nabla g(x), V(x) \rangle = 0$ having a constant sign by Proposition 4, one deduces that $\widetilde{\pi}$-almost surely, $\Omega_\infty(g)(x) = 0$. The probability $\widetilde{\pi}$ is thus carried by the equilibrium points of the dynamical system. The proposition is proved.  □

*Asymptotic independence.* In the case where $(y(t))$ has a unique equilibrium point $y_\infty$, by using the convergence of the invariant distributions $\pi_N$ to the Dirac distribution $\delta_{y_\infty}$ and the fact that the coordinates of $(X_i^N(t))$ are exchangeable, it is easy (and quite classical) to show that for any subset $I$ of coordinates, the random variables $(X_i^N(t), i \in I)$ at equilibrium are asymptotically independent with $y_\infty$ as a common limiting distribution. See Sznitman [23], for example. To summarize, the uniqueness of an equilibrium point implies that, asymptotically, the invariant distribution of the Markov process $(X_i^N(t))$ has a product form.

**5. A dimension reduction on $\mathbb{R}^K$.** In this section, a function $\phi$ on $\mathbb{R}_+^K$ is introduced such that $\rho \in \mathbb{R}_+^K$ is a zero of $\nabla \phi$ if and only if the corresponding probability distribution $\nu_\rho$ is an equilibrium point of $(y(t))$. Furthermore, it is shown that $\rho$ is a local minimum of $\phi$ if and only if $\nu_\rho$ is a local minimum of $g$ on the set of probability distributions on $\mathcal{X}$. In the next section, the function $\phi$ will be used to prove the bistability of the dynamical system $(y(t))$ in some cases.

For $\rho = (\rho_k) \in \mathbb{R}_+^K$, define

$$(19) \qquad \phi(\rho) = -\log Z(\rho) + \sum_{k=1}^{K} (\beta_k \rho_k - \alpha_k \log(\rho_k))$$

with $\alpha_k = \lambda_k/\gamma_k$, $\beta_k = (\gamma_k + \mu_k)/\gamma_k$, and $Z$ is the partition function

$$Z(\rho) = \sum_{n \in \mathcal{X}} \frac{\rho^n}{n!}.$$

PROPOSITION 6. *The probability distribution $\nu_\rho$ on $\mathcal{X}$ is an equilibrium point of the asymptotic dynamical system $(y(t))$ if and only if $\nabla \phi(\rho) = 0$.*

PROOF. Remark that, for $1 \leq k \leq K$,

$$\frac{\partial Z}{\partial \rho_k}(\rho) = \sum_{n : n + f_k \in \mathcal{X}} \frac{\rho^n}{n!},$$

so that

$$\frac{\partial \phi}{\partial \rho_k}(\rho) = \frac{\mu_k}{\gamma_k} - \frac{\lambda_k}{\rho_k \gamma_k} + \sum_{n : n + f_k \notin \mathcal{X}} \frac{\rho^n}{n!} \bigg/ \sum_{n \in \mathcal{X}} \frac{\rho^n}{n!},$$

hence this quantity is 0 if and only if the fixed point equation (10) holds. The proposition is proved. □

*Local minima of $\phi$ and $g$.* Proposition 1 has shown that an equilibrium point is necessarily a probability vector $\nu_\rho$ on $\mathcal{X}$ for some $\rho \in \mathbb{R}_+^K$. It has been shown that the function $g$ defined in Section 4 decreases along any trajectory of the dynamical system $(y(t))$ by Proposition 4 so that if it starts in the neighborhood of a local minimum of $g$, ultimately it reaches this point. At the normal scale, that is, for a finite network, it implies that, with an appropriate initial state, the state of the network $(X^N(t))$ will live for some (likely long) time in a subset of the states corresponding, up to a scaling, to this local minimum. For this reason, it is important to be able to distinguish stable from unstable equilibrium points of $(y(t))$. Due to the quite complicated expression defining $g$, it is not clear how the stability properties of the equilibrium points can be established directly with $g$. The function $\phi$ plays a key role in this respect, it reduces the complexity of the classification of the equilibrium points according to their stability properties.

Let $y \in \overset{\circ}{\mathcal{P}}(\mathcal{X})$. Taylor's formula for $g$ gives the relation, for $y' \in \mathcal{P}(\mathcal{X})$,

$$g(y') = g(y) + \langle \nabla g(y), y' - y \rangle + {}^t(y' - y)\mathcal{H}_g^y(y' - y) + o(\|y' - y\|^2),$$

where $\mathcal{H}_g^y$ is the Hessian matrix of $g$,

$$\mathcal{H}_g^y = \left(\frac{\partial^2 g}{\partial y_m \partial y_n}(y), m, n \in \mathcal{X}\right),$$

and $^t z$ is the transpose of vector $z$.

Propositions 1, 4 and 6 give the equivalence between:

— $y \in \mathcal{P}(\mathcal{X})$ is an equilibrium point;
— $y = \nu_\rho$ with $\nabla \phi(\rho) = 0$;
— $y \in \overset{\circ}{\mathcal{P}}(\mathcal{X})$ and $\langle \nabla g(y), y' - y \rangle = 0$, $\forall y' \in \mathcal{P}(\mathcal{X})$;

hence the relation

$$g(y') = g(\nu_\rho) + {}^t(y' - \nu_\rho)\mathcal{H}_g^{\nu_\rho}(y' - \nu_\rho) + o(\|y' - \nu_\rho\|^2)$$

holds.

It is assumed throughout this section that the Hessian matrix has nonzero eigenvalues at $\nu_\rho$ such that $\nabla \phi(\rho) = 0$. Consequently, for $\rho$ such that $\nabla \phi(\rho) = 0$, the probability vector $\nu_\rho$ is a local minimum of $g$, that is, a stable equilibrium point of $(y(t))$ if and only if the quadratic form associated to $\mathcal{H}_g^{\nu_\rho}$ satisfies the following property:

(20) $$^t h \mathcal{H}_g^{\nu_\rho} h \geq 0 \quad \text{for all } h = (h_n) \text{ with } \sum_{n \in \mathcal{X}} h_n = 0.$$

It will be shown in the following theorem that relation (20) is equivalent to the fact that the Hessian of $\phi$ at $\rho$ is a positive quadratic form, thereby establishing the dimension reduction for the problem of classification.

THEOREM 3 (Correspondence between the extrema of $g$ and $\phi$).

1. A vector $\rho \in \mathbb{R}_+^K$ is a local minimum of the function $\phi$ if and only if $\nu_\rho$ is a local minimum of the Lyapunov function $g$.
2. If $\rho$ is a saddle point for $\phi$, then $\nu_\rho$ is a saddle point for $g$.

PROOF. According to the above remarks, one has to study, on one hand, the sign of the quadratic form $h \to {}^t h \mathcal{H}_g^y h$ associated to $g$ at $y = \nu_\rho$, $\rho \in \mathbb{R}_+^K$, on the vector space of elements $h = (h_n) \in \mathbb{R}^{\mathcal{X}}$ such that the sum of the coordinates of $h$ is 0; and on the other hand, the sign of the quadratic form $\phi$ at $\rho$.

*The Hessian of $g$ and its quadratic form.* It is easily checked that

$$\frac{\partial^2 g}{\partial y_n \partial y_m}(y) = \frac{1}{y_n}\mathbb{1}_{\{n=m\}} - \sum_{k=1}^K n_k m_k \frac{\gamma_k}{\lambda_k + \gamma_k \langle \mathbb{I}_k, y \rangle}.$$

The quadratic form can be expressed as
$$^t h \mathcal{H}_g^y h = \sum_{n \in \mathcal{X}} \frac{h_n^2}{y_n} - \sum_{k=1}^{K} \frac{\gamma_k}{\lambda_k + \gamma_k \langle \mathbb{I}_k, y \rangle} \left( \sum_{n \in \mathcal{X}} n_k h_n \right)^2.$$

The change of variable $(h_n) \to (h_n/\sqrt{y_n})$ shows that if
$$H = \left\{ h = (h_n) \in \mathbb{R}^{\mathcal{X}} : \sum_{n \in \mathcal{X}} \sqrt{y_n} h_n = 0 \right\},$$

then it is enough to study the sign of the quadratic form $G_y$ on $H$ given by
$$G_y(h) = \sum_{n \in \mathcal{X}} h_n^2 - \sum_{k=1}^{K} \frac{\gamma_k}{\lambda_k + \gamma_k \langle \mathbb{I}_k, y \rangle} \left( \sum_{n \in \mathcal{X}} n_k \sqrt{y_n} h_n \right)^2$$
$$= \langle h, h \rangle - \sum_{k=1}^{K} \langle v_k^y, h \rangle^2,$$

where, for $1 \leq k \leq K$, $v_k^y \in R_+^{\mathcal{X}}$ is defined as
$$v_k^y = \frac{\sqrt{\gamma_k}}{\sqrt{\lambda_k + \gamma_k \langle \mathbb{I}_k, y \rangle}} (n_k \sqrt{y_n}, n \in \mathcal{X}).$$

Set
$$w_k^y \stackrel{\text{def.}}{=} \frac{\sqrt{\gamma_k}}{\sqrt{\lambda_k + \gamma_k \langle \mathbb{I}_k, y \rangle}} (\sqrt{y_n}(n_k - \langle \mathbb{I}_k, y \rangle), n \in \mathcal{X}),$$

then it is easy to check that $w_k^y$ is the orthogonal projection of $v_k^y$ in the vector space $H$, therefore
$$G_y(h) = \langle h, h \rangle - \sum_{k=1}^{K} \langle w_k^y, h \rangle^2.$$

If $W_y$ is the subvector space of $H$ generated by the vectors $w_k^y$, $1 \leq k \leq K$ and $P_{W_y}$ (resp., $P_{W_y^\perp}$) is the orthogonal projection on $W_y$ (resp., on the orthogonal of $W_y$), then

(21)
$$G_y(h) = \langle P_{W_y^\perp}(h), P_{W_y^\perp}(h) \rangle + \langle P_{W_y}(h), P_{W_y}(h) \rangle - \sum_{k=1}^{K} \langle w_k^y, P_{W_y}(h) \rangle^2$$
$$= \langle P_{W_y^\perp}(h), P_{W_y^\perp}(h) \rangle + G_y(P_{W_y}(h)).$$

To determine the sign $G_y$ on $H$, it is thus enough to have the sign of $G_y(h)$ for $h \in W_y$, such an element can be written as $h = a_1 w_1^y + \cdots + a_K w_K^y$ with $(a_k) \in \mathbb{R}^K$,
$$G_y(h) = \sum_{1 \leq i,j \leq K} a_i a_j \left( \langle w_i^y, w_j^y \rangle - \sum_{k=1}^{K} \langle w_k^y, w_i^y \rangle \langle w_k^y, w_j^y \rangle \right)$$

hence, if $\mathcal{W}^y$ is the $K \times K$ matrix defined by $\mathcal{W}^y = (\langle w_k^y, w_l^y \rangle, 1 \leq k, l \leq K)$,

$$(22) \qquad G_y(h) = {}^t a \mathcal{W}^y (I - \mathcal{W}^y) a \qquad \text{for } h = \sum_{k=1}^{K} a_k w_k^y.$$

The eigenvalues of the matrix $\mathcal{W}^y$ being all real (it is symmetrical) and nonnegative since its associated quadratic form is nonnegative, therefore $G_y$ is positive on $W_y$ if and only if all the eigenvalues of $\mathcal{W}^y$ are in the interval $(0, 1)$.

*The Hessian of $\phi$ and its quadratic form.* For $\rho \in \mathbb{R}_+^K$,

$$\frac{\partial^2 \phi}{\partial \rho_k \, \partial \rho_l}(\rho) = -\frac{\partial^2 \log Z}{\partial \rho_k \, \partial \rho_l}(\rho) + \frac{\alpha_k}{\rho_k^2} \mathbb{1}_{\{k=l\}}$$

with $(\alpha_k) = (\lambda_k/\gamma_k)$, for $1 \leq k, l \leq K$. The derivatives of $\log Z$ have the following properties:

$$(23) \qquad \rho_k \frac{\partial \log Z}{\partial \rho_k}(\rho) = \frac{1}{Z(\rho)} \sum_{n \in \mathcal{X}} n_k \frac{\rho^n}{n!} = \langle \mathbb{I}_k, \nu_\rho \rangle$$

and

$$-\rho_k \rho_l \frac{\partial^2 \log Z}{\partial \rho_k \, \partial \rho_l}(\rho) = -\frac{1}{Z(\rho)} \sum_{n \in \mathcal{X}} n_k n_l \frac{\rho^n}{n!}$$

$$+ \left( \frac{1}{Z(\rho)} \sum_{n \in \mathcal{X}} n_k \frac{\rho^n}{n!} \right) \left( \frac{1}{Z(\rho)} \sum_{n \in \mathcal{X}} n_l \frac{\rho^n}{n!} \right)$$

$$+ \mathbb{1}_{\{k=l\}} \frac{1}{Z(\rho)} \sum_{n \in \mathcal{X}} n_k \frac{\rho^n}{n!},$$

hence

$$-\rho_k \rho_l \frac{\partial^2 \log Z}{\partial \rho_k \, \partial \rho_l}(\rho) = \langle \mathbb{I}_k, \nu_\rho \rangle \langle \mathbb{I}_l, \nu_\rho \rangle - \langle \mathbb{I}_{k,l}, \nu_\rho \rangle + \mathbb{1}_{\{k=l\}} \langle \mathbb{I}_k, \nu_\rho \rangle,$$

where

$$\langle \mathbb{I}_{k,l}, \nu_\rho \rangle = \frac{1}{Z(\rho)} \sum_{n \in \mathcal{X}} n_k n_l \frac{\rho^n}{n!}.$$

The quadratic form associated to $\phi$ at $\rho \in \mathbb{R}_+^K$ is given by, for $a = (a_k) \in \mathbb{R}^K$,

$$\Phi_\rho(a) = \sum_{1 \leq k, l \leq K} (\langle \mathbb{I}_k, \nu_\rho \rangle \langle \mathbb{I}_l, \nu_\rho \rangle - \langle \mathbb{I}_{k,l}, \nu_\rho \rangle) \frac{a_k \, a_l}{\rho_k \, \rho_l} + \sum_{k=1}^{K} (\alpha_k + \langle \mathbb{I}_k, \nu_\rho \rangle) \frac{a_k^2}{\rho_k^2}.$$

By using the change of variable (recall that $\alpha_k = \lambda_k/\gamma_k$),

$$a = (a_k) \to \left(\frac{\sqrt{\lambda_k + \gamma_k \langle \mathbb{I}_k, \nu_\rho \rangle}}{\sqrt{\gamma_k}} \frac{a_k}{\rho_k}\right),$$

one gets that the sign of $\phi_\rho$ has the same range as the sign of $\Psi_\rho$, where

$$\Psi_\rho(a) = \langle a, a \rangle + \sum_{1 \leq k,l \leq K} \frac{\sqrt{\gamma_k}}{\sqrt{\lambda_k + \gamma_k \langle \mathbb{I}_k, \nu_\rho \rangle}} \frac{\sqrt{\gamma_l}}{\sqrt{\lambda_l + \gamma_l \langle \mathbb{I}_l, \nu_\rho \rangle}}$$

$$\times (\langle \mathbb{I}_k, \nu_\rho \rangle \langle \mathbb{I}_l, \nu_\rho \rangle - \langle \mathbb{I}_{k,l}, \nu_\rho \rangle) a_k a_l$$

$$= \langle a, a \rangle - \sum_{1 \leq k,l \leq K} \langle w_k^{\nu_\rho}, w_l^{\nu_\rho} \rangle a_k a_l,$$

with the above notation. Therefore, the sign of the quadratic form associated to $\phi$ at $\rho$ has the same values as the sign of $\Psi_\rho(a)$ defined by

(24) $$\Psi_\rho(a) = {}^t a (I - \mathcal{W}^{\nu_\rho}) a, \qquad a = (a_k) \in \mathbb{R}^K.$$

Equations (22) and (24) show that $G_{\nu_\rho}$ is positive on $W_{\nu_\rho}$ if and only if $\Psi_\rho$ is positive on $\mathbb{R}_+^K$ which proves assertion 1 of the theorem. Similarly, if $\rho$ is a saddle point of $\phi$, equation (24) shows that the matrix $\mathcal{W}^{\nu_\rho}$ has eigenvalues in $(0,1)$ and in $(1,+\infty)$, so that $G_{\nu_\rho}$ takes positive and negative values on $W$, and hence on $H$, $\nu_\rho$ is thus a saddle point of $g$. The theorem is proved. □

**6. Bistability of the asymptotic dynamical system.** This section gives an example where the asymptotic dynamical system has at least three fixed points: Two of them are stable and the other is a saddle point. The corresponding stochastic network therefore exhibits a metastability property. In the limit, it suggests that its state switches from one stable point to the other after a long residence time. See Figure 2. The problem of estimating the residence time in the neighborhood of a stable point is not addressed here. According to examples from statistical physics, the expected value of this residence time should be of exponential order with respect to the size $N$ of the network. For reversible Markov processes, Bovier [3, 4, 5] presents a potential theoretical approach to get lower and upper bounds for this expected value.

*A network with two classes.* A simple setting is considered here: There are two classes of customers, $K = 2$, the capacity requirements are $A_1 = 1$ (small customers) and $A_2 = C$ (large ones) so that, at a given node, there may be $n$ class 1 customers, $0 \leq n \leq C$, or only one class 2 customer. It is assumed that $\gamma_1 = \gamma_2 = 1$ and $\mu_1 = \mu_2 = 0$ so that a customer leaves the network only when it is rejected at some node.

The two classes cannot coexist at a given node and, moreover, when a node contains class 1 customers, it has to get completely empty before accommodating a class 2 customer. Moreover, when the network is mostly filled with class 1 customers, the competition for capacity at each node should be favorable to class 1 customers, due to their large internal arrival rate (i.e., their arrival rate from all the other nodes). This can explain the stability of a state with a high density in class 1 customers. The same intuitive argument holds for the existence of a stable state with a comparatively higher density in class 2 customers, though it is clear that the occurrence of this phenomenon should depend on the compared values of the different arrivals, services and transfers rates.

PROPOSITION 7. *For a network with two classes of customers such that $A_1 = 1$, $A_2 = C$, $\gamma_1 = \gamma_2 = 1$, $\mu_1 = \mu_2 = 0$, for $C$ sufficiently large, there exist $\lambda_1$ and $\lambda_2 \in \mathbb{R}_+$ such that the corresponding energy function $\phi$ has at least one saddle point and two local minima.*

From Theorem 3, one deduces that there exists a stochastic network whose asymptotic dynamical system has at least two stable points.

PROOF OF PROPOSITION 7. Fix $\rho \in \mathbb{R}_+^2$ and choose $(\lambda_1, \lambda_2) \in \mathbb{R}_+^2$ so that $\rho$ satisfies equations (8), that is,

$$(25) \qquad \lambda_k = \rho_k - \langle \mathbb{I}_k, \nu_\rho \rangle = \rho_k \left( 1 - \frac{\partial \log Z}{\partial \rho_k}(\rho) \right), \qquad k = 1, 2,$$

by relation (23), so that $\nu_\rho$ is an equilibrium point for the limiting dynamics. It will be assumed for the moment that $C = +\infty$. The corresponding function $\phi$ is then given by

$$\tilde{\phi}(\rho) = -\log(\rho_2 + e^{\rho_1}) + \rho_1 + \rho_2 - \lambda_1 \log \rho_1 - \lambda_2 \log \rho_2.$$

Using equation (25), one gets that

$$\frac{\partial^2 \tilde{\phi}}{\partial \rho_1^2}(\rho) = \frac{\lambda_1}{\rho_1^2} - \frac{\rho_2 e^{\rho_1}}{(\rho_2 + e^{\rho_1})^2}$$

$$= \frac{\rho_2(\rho_2 + (1 - \rho_1)e^{\rho_1})}{\rho_1(\rho_2 + e^{\rho_1})^2}$$

and

$$\frac{\partial^2 \tilde{\phi}}{\partial \rho_2^2}(\rho) = \frac{\lambda_2}{\rho_2^2} + \frac{1}{(\rho_2 + e^{\rho_1})^2} > 0.$$

If $\bar{\rho} = (\bar{\rho}_1, \bar{\rho}_1)$ is chosen such that the inequality $\bar{\rho}_2 < (\bar{\rho}_1 - 1) \exp(\bar{\rho}_1)$ holds, then

$$\frac{\partial^2 \tilde{\phi}}{\partial \rho_1^2}(\bar{\rho}) < 0 \quad \text{and} \quad \frac{\partial^2 \tilde{\phi}}{\partial \rho_2^2}(\bar{\rho}) > 0.$$

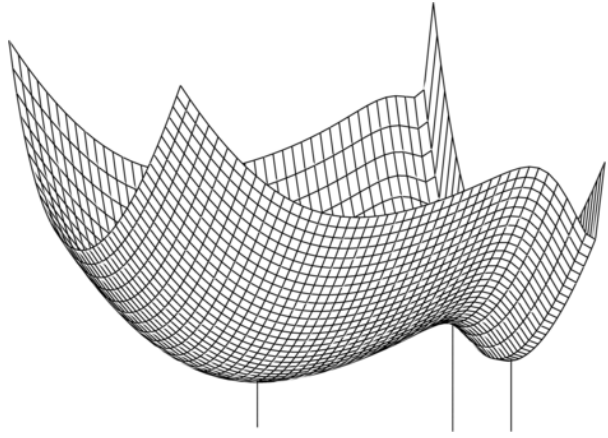

FIG. 1. *Function $\phi$ with one saddle point and two stable equilibrium points. Two classes with $\lambda_1 = 0.68$, $\lambda_2 = 9.0$, $A_1 = 1$ and $A_2 = C = 20$.*

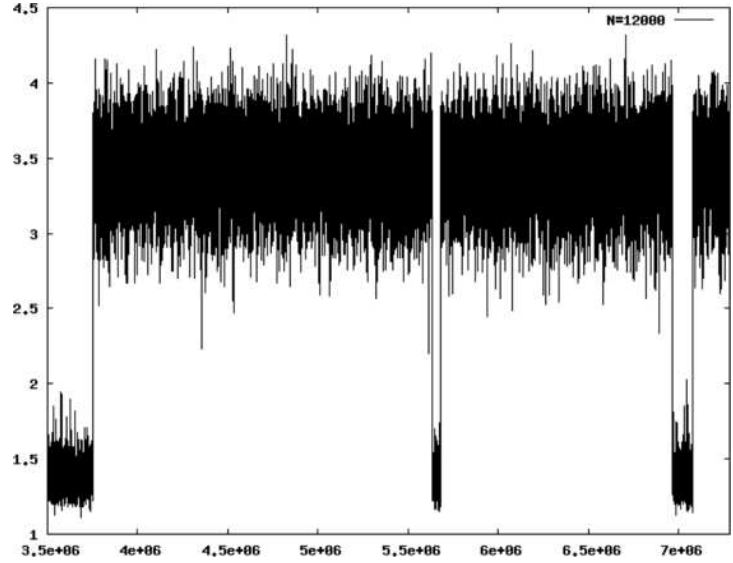

FIG. 2. *Time evolution of the proportion of nodes without class 2 particle. Case $N = 12000$ nodes, $A_1 = 1$, $A_2 = C = 5$, $\lambda_1 = 0.64$, $\lambda_2 = 2.71$, $\mu_1 = \mu_2 = 0$ and $\gamma_1 = \gamma_2 = 1$.*

The constant $C$ is now assumed to be finite and sufficiently large so that the above inequalities with $\phi$ in place of $\tilde{\phi}$ are satisfied, $\bar{\rho}$ is a saddle point for $\phi$. The function $\phi$ is given by

$$\phi(\rho) = -\log\left(\rho_2 + \sum_{n=0}^{C} \frac{\rho_1^n}{n!}\right) + \rho_1 + \rho_2 - \lambda_1 \log \rho_1 - \lambda_2 \log \rho_2.$$

The function $\rho_2 \to \phi(\bar{\rho}_1, \rho_2)$ is convex, $\bar{\rho}_2$ is a strict local minimum by construction and therefore a *global* minimum. Similarly, the function $\rho_1 \to \phi(\rho_1, \bar{\rho}_2)$ has a strict local maximum at $\bar{\rho}_1$,

$$\inf\{\phi(\rho): \rho = (\rho_1, \bar{\rho}_2), \rho_1 < \bar{\rho}_1\} < \phi(\bar{\rho}),$$
$$\inf\{\phi(\rho): \rho = (\rho_1, \bar{\rho}_2), \bar{\rho}_1 < \rho_1\} < \phi(\bar{\rho}) = \inf\{\phi(\rho): \rho \in \Delta\},$$

with $\Delta = \{(\bar{\rho}_1, \rho_2): \rho_2 \in \mathbb{R}_+ \setminus \{0\}\}$. Since $\phi((\rho_1, \rho_2))$ converges to $+\infty$ when $\rho_1$ or $\rho_2$ converges to 0 or $+\infty$, one concludes that the function $\phi$ has at least two local finite minima, one on each side of $\Delta$. The proposition is proved. Figure 1 gives an example of such a situation. $\square$

N. ANTUNES
FACULDADE DE CIÊNCIAS E TECNOLOGIA
UNIVERSIDADE DO ALGARVE
CAMPUS DE GAMBELAS
8005-139 FARO
PORTUGAL
E-MAIL: nantunes@ualg.pt

C. FRICKER
PH. ROBERT
INRIA
DOMAINE DE VOLUCEAU
B.P. 105
78153 LE CHESNAY CEDEX
FRANCE
E-MAIL: Christine.Fricker@inria.fr
Philippe.Robert@inria.fr

D. TIBI
UNIVERSITÉ PARIS 7
UMR 7599
2 PLACE JUSSIEU
75251 PARIS CEDEX 05
FRANCE
E-MAIL: Danielle.Tibi@math.jussieu.fr